\documentclass{article}
\usepackage{latexsym,amssymb,amsmath,amsthm,amsfonts,euscript}

\newcounter{theorem}

\newcounter{proposition}

\begin{document}

\begin{center}
{\Large \bf An analog of the Hille theorem for hypercomplex functions in a finite-dimensional commutative algebra}

\vskip 2mm

{\bf \large S.A. Plaksa, V.S. Shpakivskyi and M.V. Tkachuk}

\vskip 2mm

{Institute of Mathematics of NAS of Ukraine, Kyiv, Ukraine}

\vskip 2mm

{plaksa62@gmail.com}\\
shpakivskyi86@gmail.com\\
maxim.v.tkachuk@gmail.com
\end{center}

\begin{abstract}
\noindent 
We prove that a locally bounded and differentiable in the sense of G\^ateaux function given in a finite-dimensional commutative Banach algebra over the complex field is also differentiable in the sense of Lorch.
\end{abstract}





\section{Introduction}

For mappings of vector spaces as well as commutative Banach algebras, various concepts of differentiable mappings are used. Let us consider
some of these concepts and the relationship between them (for more detailed information, see, for example, the monograph \cite{Pla-Shpak-mono}).

Let $V_j$, $j=1,2$, be normalized vector spaces with the norms $\|\cdot\|_{V_j}$, and let $\Omega$ be an open subset of $V_1$.

A mapping $\Phi \colon \Omega \rightarrow V_2$ is
{\em Fr\'echet differentiable}  at a point $\zeta\in\Omega$ if there exists a bounded
linear operator $A_{\zeta} \colon V_1\rightarrow V_2$ such that
\[\frac{\left\|\Phi(\zeta+h)-\Phi(\zeta)-A_{\zeta}\,h\right\|_{V_2}}{\|h\|_{V_1}}\to 0\,,\quad \|h\|_{V_1}\to 0\,.\]
The operator\, $A_{\zeta}$\, is called the {\em Fr\'echet derivative} of the mapping $\Phi$ at the point $\zeta$  (cf. M.~Fr\'echet \cite{Freche-1-1}).

 A mapping $\Phi \colon \Omega \rightarrow V_2$ is
{\em  G\^ateaux differentiable} at a point $\zeta\in\Omega$ if the {\em G\^ateaux differential}
\begin{equation} \label{1-1:Gdif}
\mathcal{D}_G\Phi(\zeta,h):=\lim\limits_{\delta\rightarrow0}\frac{\Phi(\zeta+\delta h)-\Phi(\zeta)}{\delta}\,,
\end{equation}
exists for all $h\in V_1$, where $\delta$ is taken from the scalar field associated with the space $V_1$ (cf. R.~G\^ateaux \cite{Gato-1-1}).
In addition,
if there exists a bounded
linear operator $B_{\zeta} \colon V_1\rightarrow V_2$ such that
$\mathcal{D}_G\Phi(\zeta,h)\equiv B_{\zeta}\,h$\,, the operator\,
$B_{\zeta}$\, is called the {\em G\^ateaux derivative}
of the mapping $\Phi$ at the point $\zeta$.

It is evident that if a mapping $\Phi \colon \Omega \rightarrow V_2$ is Fr\'echet differentiable at a point
 $\zeta\in\Omega$, then the G\^ateaux derivative of $\Phi$ exists at the point $\zeta$ and is equal to the
 Fr\'echet derivative of the mapping $\Phi$ at that point.
The converse is not true in general.

The following theorem, which is proved by E.~Hille \cite{Hille-44} (see also E.~Hille and R.~Phillips \cite[p.~112]{Hil_Filips-1-1}), indicates conditions sufficient for the G\^ateaux derivative to coincide with the Fr\'echet derivative.
\vskip 1mm

\textbf{Hille theorem.} {\it
Let $V_1$, $V_2$ be complex Banach spaces and $\Omega$ be an open subset of\, $V_1$. Suppose that
a mapping $\Phi \colon \Omega \rightarrow V_2$ is locally bounded in $\Omega$ and
has the G\^ateaux derivative $B_{\zeta}$ at every point $\zeta\in\Omega$. Then $B_{\zeta}$ is the Fr\'echet derivative of the mapping $\Phi$
at every point $\zeta\in\Omega$.}
\vskip 1mm

The considered Fr\'echet and G\^ateaux derivatives are defined as bounded linear operators.
Considering mappings given in commutative Banach algebras, it is possible to introduce
the concepts of derivatives that are understood as functions defined in the same domain as a
given function.   

First of all, note that for functions $\Phi(\zeta)$ given in a domain of a finite-dimensional
algebra, G.~Scheffers \cite{Scheffers-1-2} considered a derivative $\Phi'(\zeta)$ defined by $d\Phi=\Phi'(\zeta)d\zeta$,
which is understood as a hypercomplex function.

Generalizing such an approach to the case of functions given in a
domain of an arbitrary commutative Banach algebra,
E.R.~Lorch \cite{Lorch-1-2} introduced a derivative, which is also
understood as a function given in the same domain.

Let $\mathbb A$ be a
commutative Banach algebra with unit $1$ over either the field of real numbers $\mathbb R$ or the field of
complex numbers $\mathbb C$\,, and let $\|\zeta\|$ denote the norm of element $\zeta\in\mathbb A$.

Let us fix a finite number $k$ of the vectors $e_1=1,e_2,\dots,e_k\in\mathbb{A}$ so that they were linearly independent over the field of real numbers
$\mathbb{R}$.
Let
\begin{equation}\label{intr-E_k}
E_k:=\{\zeta=\sum\limits_{j=1}^k x_{j}e_j:\,\,\,x_j\in\mathbb{R}\}
\end{equation}
be the linear span of the vectors\,\,
$e_1,e_2,\dots,e_k$\,\, over the field $\mathbb{R}$\,.

\vskip 2mm

\textbf{Definition 1.}
A function $\Phi \colon \Omega \rightarrow \mathbb{A}$ given in a domain\, $\Omega\subset E_k$\,
is called {\em differentiable in the sense of Lorch at a point}\, $\zeta\in\Omega$\,
if there exists an element
$\Phi_L'(\zeta)\in\mathbb{A}$ such that for each $\varepsilon>0$
there exists $\delta>0$ such that for all\, $h\in E_k$\, with\, $\|h\|
<\delta$\, the following inequality is fulfilled:
\[\left\|\Phi(\zeta+h)-\Phi(\zeta)-h\Phi_L'(\zeta)\right\|\leq \|h\|\,\varepsilon\,.\]
Here, $\Phi_L'(\zeta)$\, is called the {\em Lorch derivative} of the function $\Phi$ at the point $\zeta$  (cf. E.R.~Lorch
\cite{Lorch-1-2}).

\vskip 2mm

It is clear that a function $\Phi$, which is differentiable in the sense
of Lorch at a point\, $\zeta\in\Omega$\,,
is also the Fr\'echet differentiable at the same point. The converse is not true.

Using the G\^ateaux differential \eqref{1-1:Gdif},  I.P.~Mel'nichenko
\cite{Mel'nichenko75-1-2} defined a derivative that is also a hypercomplex function.
\vskip 2mm

\textbf{Definition 2.}
A function $\Phi \colon \Omega \rightarrow \mathbb{A}$ given in a
domain\, $\Omega\subset E_k$\, is called {\em differentiable in the
sense of G\^ateaux at a point}\,  $\zeta\in\Omega$\, if there exists an
element $\Phi_G'(\zeta)\in\mathbb{A}$ such that
\[\lim\limits_{\delta\rightarrow 0^+} \frac{\Phi(\zeta+\delta
h)-\Phi(\zeta)}{\delta}= h\,\Phi_G'(\zeta)\qquad\forall\,h\in
E_k.\]
We call\, $\Phi_G'(\zeta)$\, the {\em G\^ateaux--Mel'nichenko derivative}
of the function $\Phi$ at the point $\zeta$ (cf. I.P.~Mel'nichenko \cite{Mel'nichenko75-1-2}).

\vskip 2mm

It is clear that if a function $\Phi$ is differentiable in the sense of
Lorch in a domain\, $\Omega\subset E_k$\,, then it is also differentiable in the sense of
G\^ateaux, and $\Phi_L'(\zeta)=\Phi_G'(\zeta)$ for all
$\zeta\in\Omega$. The converse is clearly not true similarly to
the fact that the existence of all di\-rec\-ti\-o\-nal derivatives at a
point does not guarantee a strong differentiability (or even
continuity) of function at that point.
Moreover, the Fr\'echet differentiability of
a function $\Phi \colon \Omega \rightarrow \mathbb{A}$ does not imply the differentiability of this function
in the sense of G\^ateaux (see Examples 3.1 and 3.2 in \cite{Pla-Shpak-mono}).

In the paper \cite{Pl-zb17-2-3},
an analog of the Hille theorem is proved for
hypercomplex functions in a concrete three-dimensional
commutative algebra over the complex field (see also Theorem 6.17 in \cite{Pla-Shpak-mono}).
Namely, it is proved that locally bounded and
differentiable in the sense of G\^ateaux functions given in a three-dimensional commutative
harmonic algebra with two-dimensional radical are also differentiable in the sense of Lorch.
Note that it is impossible to establish this result using the Hille theorem because
the Fr\'echet differentiability does not imply the existence of Lorch derivative.

The purpose of this paper is to prove a similar theorem in an arbitrary finite-dimensional commutative Banach algebra.

\vskip 2mm

\section{A finite-dimensional commutative associative algebra and the Cartan basis}

E.~Cartan \cite{Cartan-3-1} proved that for an arbitrary $n$-dimensional
commutative associative algebra $\mathbb{A}$ with unit over the field of complex number $\mathbb{C}$,
  there exists a basis  $\{I_k\}_{k=1}^{n}$ and there
  exist structural constants $\Upsilon_{r,k}^{s}$ such that
the following multiplication rules hold:
\begin{align}
1.  & \hspace*{2mm}  \forall\, r,s\in[1,m]\cap\mathbb{N}\,: &
I_rI_s &=\left\{
\begin{array}{rcl}
0 &\hspace*{2mm} \mbox{for} & \hspace*{2mm} r\neq s,\\[2mm]
I_r &\hspace*{2mm} \mbox{for} & \hspace*{2mm} r=s;\\
\end{array}
\right. \notag \\[5mm]
2. & \hspace*{2mm} \forall\, r,s\in[m+1,n]\cap\mathbb{N}\,: &
I_rI_s &=
\sum\limits_{k=\max\{r,s\}+1}^n\Upsilon_{r,k}^{s}I_k\,;\notag
\end{align}
\begin{multline*} 
3. \hspace*{2mm} \forall\,s\in[m+1,n]\cap\mathbb{N}\;  \hspace{2mm} \exists!\;
 u_s\in[1,m]\cap\mathbb{N}\;  \hspace*{2mm} \forall\, r\in[1,m]\cap\mathbb{N}\,:\\
I_rI_s=\left\{
\begin{array}{ccl}
0 & \hspace*{2mm} \mbox{for} & \hspace*{2mm} r\neq u_s\,,\\[2mm]
I_s & \hspace*{2mm} \mbox{for} & \hspace*{2mm} r= u_s\,, \\
\end{array}
\right.
\end{multline*}
where $\mathbb{N}$ is the set of natural numbers.

%


It is obvious that the first $m$ basic vectors $I_1,I_2,\dots,I_m$ are idempotents,
and the vectors
$I_{m+1},I_{m+2},\dots,I_n$ 
are nilpotent elements.

In what follows, the algebra $\mathbb{A}$ with Cartan basis is
denoted as $\mathbb{A}_n^m$.
 The element $1=I_1+I_2+\dots+I_m$
is the unit in the algebra\, $\mathbb{A}_n^m$\, that is a commutative Banach algebra with
Euclidean norm defined by the equality
$$\|v\|:=\sqrt{\sum\limits_{j=1}^n|v_j|^2}\,,\qquad \mbox{where}\quad v=\sum\limits_{j=1}^n v_jI_j, \quad v_i\in\mathbb{C}\,.$$

The algebra $\mathbb{A}_n^m$ contains $m$ maximal ideals
$$\mathcal{I}_u:=\Biggr\{\sum\limits_{r=1,\,r\neq u}^n\lambda_rI_r:\lambda_r\in
\mathbb{C}\Biggr\}, \quad  u=1,2,\ldots,m,
$$
and their intersection is the radical $$\mathcal{R}:=
\Bigr\{\sum\limits_{r=m+1}^n\lambda_rI_r:\lambda_r\in
\mathbb{C}\Bigr\}.$$

Consider $m$ linear functionals
$f_u:\mathbb{A}_n^m\rightarrow\mathbb{C}$ satisfying the
equalities
\[f_u(I_u)=1,\quad f_u(\omega)=0\quad\forall\,\omega\in\mathcal{I}_u\,,
\quad u=1,2,\ldots,m.\]
Inasmuch as the kernel of functional $f_u$ is the maximal ideal
$\mathcal{I}_u$, this functional is also continuous and
multiplicative (see \cite[p. 147]{Hil_Filips-1-1}).

\vskip 2mm

\section{The main result}

Consider vectors $e_1=1,e_2,\ldots,e_k\in\mathbb{A}_n^m$,\,\, $2\leq k\leq2n$, which are linearly independent over the field 
$\mathbb{R}$. It means that the equality
$$\sum\limits_{j=1}^k\alpha_je_j=0,\qquad \alpha_j\in\mathbb{R},$$
holds if and only if $\alpha_j=0$ for all $j=1,2,\ldots,k$.
The following decompositions with respect to
the Cartan basis $\{I_r\}_{r=1}^n$ hold:
\[e_1=\sum\limits_{r=1}^mI_r\,, 
\quad e_j=\sum\limits_{r=1}^na_{jr}\,I_r\,,\quad a_{jr}\in\mathbb{C},\quad j=2,3,\ldots,k.\]


The linear span $E_k$ is defined by equality \eqref{intr-E_k}.
In what follows, we impose the following restriction on the choice of the linear span\, $E_k$:
\begin{equation}\label{cond-on-Ek}
 \{f_u(\zeta) : \zeta\in E_k\} =\mathbb{C}\,,\qquad   u=1,2,\ldots,m\,,
\end{equation}
i.e., the images of the set\, $E_k$\, under all mappings\, $f_u$\, must be the whole complex plane (cf. \cite{Pl-Pukh-Analele}).
 Obviously, it holds if and only if for every fixed $u=1,2, \ldots, m$
at least one of the numbers $a_{2u}$, $a_{3u},\ldots,a_{ku}$ belongs to
$\mathbb{C}\setminus\mathbb{R}$.


\vskip 2mm

\textbf{Definition 3.}
We say that a function $\Phi \colon \Omega  \longrightarrow
\mathbb{A}_n^m$ is \textit{monogenic} in a domain
$\Omega\subset E_k$ if $\Phi$ is continuous and
differentiable in the sense of G\^ateaux at every point of
$\Omega$.

\vskip 2mm


In \cite{Pla-Shpak-mono,Sh-co} for the case $k=3$, we developed a theory of monogenic functions $\Phi \colon \Omega  \longrightarrow
\mathbb{A}_n^m$, which includes analogs of classical theorems of complex analysis (the Cauchy integral theorems for a curvilinear and for a surface integral, the Cauchy integral formula, the Morera theorem, the Taylor theorem).
These results are generalized in the papers \cite{Sh-co-Zb,Shpakivskyi-Zb-2015} for the case $2\leq k\leq2n$.

Note that we use the notion of monogenic function in the sense of existence
of derived numbers for this function (cf. the monographs by E.~Goursat
\cite{Goursat-1-2} and Ju.Ju.~Trokhimchuk \cite{Trokhimchuk-1-2} ) in a combination with its continuity.
In the scientific literature the denomination of monogenic
function is used else for functions given in non-commutative
algebras and satisfying certain conditions similar to the
classical Cauchy--Riemann conditions (see, for example,
F.~Sommen \cite{Sommen-81} and J.~Ryan \cite{Ryan-1-2}).


The main result of this paper is the following

\vskip 2mm

\textbf{Theorem 1 (an analog of the Hille theorem).} \label{teo_1_konstruct_opys_A_n_m-k}
 \textit{ Suppose that condition \eqref{cond-on-Ek} is satisfied.
For a function $\Phi:\Omega\rightarrow\mathbb{A}_n^m$ given in an arbitrary domain $\Omega \subset
E_k$ the following properties are equivalent:
\begin{description}
  \item[(I)] $\Phi$ is a function
locally bounded and differentiable in the sense of G\^ateaux
in $\Omega $;
  \item[(II)] $\Phi$ is a monogenic function in $\Omega $;
  \item[(III)]  $\Phi$ is a function differentiable in the sense of Lorch
in $\Omega $.
\end{description}
}

\vskip 2mm

It follows from Theorem 1 
that under condition \eqref{cond-on-Ek} the property of continuity in Definition 3 of a monogenic function can be replaced by its local boundedness. Note that in Theorem 2 in \cite{Tk-Pl}, which is an analog of Menchov--Trokhimchuk theorem,
the conditions of monogeneity are weakened in an other way for continuous functions taking values in one of
three-dimensional commutative algebras over the complex field (cf. also Theorem 6.18 in \cite{Pla-Shpak-mono}).

Certainly, Theorem 1 
yields that the property of function to be locally bounded and
differentiable in the sense of G\^ateaux in $\Omega$ is also
equivalent to the various definitions of monogenic function, which
are stated in Theorem 9.9 in \cite{Pla-Shpak-mono}.

\vskip 2mm

\section{Proof of the main result}

\subsection{Auxiliary statements}

In what follows,
\begin{equation}\label{zeta}
\zeta:=\sum\limits_{j=1}^kx_j\,e_j\,,\qquad  x_j\in\mathbb{R}.
\end{equation}
It is obvious that
\begin{equation}\label{xi-u}
 \xi_u:=f_u(\zeta)=x_1+\sum\limits_{j=2}^kx_j\,a_{ju},\qquad u=1,2,\ldots,m.
\end{equation}

The following expansion of resolvent is proved in Lemma 1 \cite{Sh-co-Zb}:
\begin{equation}\label{lem-rez}
(te_1-\zeta)^{-1}=\sum\limits_{u=1}^m\frac{1}{t-\xi_u}\,I_u+
 \sum\limits_{s=m+1}^{n}\sum\limits_{r=2}^{s-m+1}\frac{Q_{r,s}}
 {\left(t-\xi_{u_{s}}\right)^r}\,I_{s}\,
 \end{equation}
$$\forall\,t\in\mathbb{C}:\,
t\neq \xi_u,\quad u=1,2,\ldots,m,$$
where the coefficients\ $Q_{r,s}$ are determined by the following recurrence
relations:
\[\begin{array}{c}
\displaystyle
Q_{2,s}=T_s\,,\quad
Q_{r,s}=\sum\limits_{q=r+m-2}^{s-1}Q_{r-1,q}\,B_{q,\,s}\,,\; \;\;r=3,4,\ldots,s-m+1,\\
\end{array} \]
and
\[T_s:=\sum\limits_{j=2}^kx_ja_{js}\,, \quad s=m+1,m+3,2\ldots,n,\]
\[B_{q,s}:=\sum\limits_{p=m+1}^{s-1}T_p \Upsilon_{q,s}^p\,,
\quad s=m+2,m+3,\ldots,n,\]
 and the natural numbers $u_s$ are defined in the
rule  3  of the multiplication table of algebra
$\mathbb{A}_n^m$. Expansion \eqref{lem-rez} generalizes the corresponding expansion obtained in \cite{Sh-co}
(see also Section 8.3 in \cite{Pla-Shpak-mono}) for the case $k=3$.

It follows from expansion (\ref{lem-rez}) that
\[
\zeta^{-1}=\sum\limits_{u=1}^m\frac{1}{\xi_u}\,I_u-
 \sum\limits_{s=m+1}^{n}\sum\limits_{r=2}^{s-m+1}(-1)^r\frac{Q_{r,s}}
 {\left(\xi_{u_{s}}\right)^r}\,I_{s}\,,
\]
and the union $\bigcup\limits_{1 \leq u \leq m} M_u$ of the sets $M_u := \{ \zeta\in E_k : \xi_u = 0 \}$
is the set of all noninvertible elements in $E_k$.

Since the equality (\ref{zeta}) establishes a one-to-one correspondence between the points $\zeta\in E_k$ and $(x_1,x_2,\ldots,x_k)\in\mathbb{R}^k$, the linear span (\ref{intr-E_k}) is in fact a $k$-dimensional real vector space. Therefore, we shall use some geometric concepts (parallelism, convexity, $(k-2)$-plane etc.) of real vector spaces for
the objects from $E_k$.


Taking into account the equality (\ref{xi-u}), we can state that
the set $M_u$ of noninvertible elements \eqref{zeta}
is $(k-2)$-dimensional linear subspace of $E_k$, which is
defined by the equalities
  \[ 
x_1+\sum\limits_{j=2}^kx_j\,{\rm Re}\,a_{ju}=0,\qquad
\sum\limits_{j=2}^kx_j\,{\rm Im}\,a_{ju}=0 \]
for $u=1,2,\ldots,m$.

We shall say that $(k-2)$-plane $L\subset E_k$ is parallel to $M_u$ if $L$ can be obtained by translation of $M_u$.

Denote the line segment with the starting point $\zeta_{1}\in E_k$ and the end point $\zeta_{2}\in E_k$ as $s[\zeta_{1}, \zeta_{2}]$.
We shall use the same denotation of the line segment $s[\xi_{1}, \xi_{2}]$ in the case where $\xi_1, \xi_2\in\mathbb{C}$.

Consider the expansion
\begin{equation} \label{Phi-expansion}
\Phi(\zeta)= \sum\limits_{r=1}^{n} V_r(\zeta) I_r
\end{equation}
of a function $\Phi \colon\Omega\rightarrow\mathbb{A}_n^m$ with respect to the basis $\{I_r\}_{r=1}^{n}$, where $V_r : \Omega\rightarrow\mathbb{C}$.

 \vskip2mm

\textbf{Lemma 1.}\label{lem-0}
\textit{Let a function $\Phi \colon\Omega\rightarrow\mathbb{A}_n^m$ be differentiable in the sense of
G\^ateaux at a point $\zeta_0$ of a domain $\Omega\subset E_k$.}

(a) \textit{For any $u\in\{1,2,\dots,m\}$,
let $\Sigma_u\subset\Omega$ be a two-dimensional surface  containing four segments
$s[\zeta_{0}, \zeta_{1}]$, $s[\zeta_{0}, \zeta_{2}]$, $s[\zeta_{0}, \zeta_{3}]$, $s[\zeta_{0}, \zeta_{4}]$ such that their images
under the mapping $f_{u}$ are subsets of the segments $s[\xi_0,\xi_0 + 1]$, $s[\xi_0,\xi_0 + i]$, $s[\xi_0,\xi_0 - 1]$, $s[\xi_0,\xi_0 - i]$, respectively, where $\xi_0 = f_{u}(\zeta_0)$, and, in addition, let the restriction of the functional $f_{u}$ to the surface $\Sigma_u$ carry out a one-to-one mapping of this surface onto the domain $Q_u:=\{\xi=f_{u}(\zeta) : \zeta\in \Sigma_u\}$.
Then the function $H_u \colon Q_u \rightarrow \mathbb{C}$ defined as $H_u(\xi) := V_u(\zeta)$ for all
$\xi = \tau + i \eta=f_{u}(\zeta)$ for $\zeta\in \Sigma_u$ and $\tau,\eta\in\mathbb{R}$,
where the function $V_u$ is defined in expansion \eqref{Phi-expansion},  has the partial derivatives $\partial H_u/\partial \tau$, $\partial H_u/\partial \eta$ at the point $\xi_0$, and the following equality holds:}
\begin{equation} \label{clas-CR-cond}
\frac{\partial H_u}{\partial \eta} = i \frac{\partial H_u}{\partial \tau}.
\end{equation}

(b) \textit{In the case where the function $\Phi$ takes values in the radical $\mathcal{R}$, i.e., expansion \eqref{Phi-expansion} takes the form
\begin{equation}\label{Phi-expansion-rad}
\Phi(\zeta)= \sum\limits_{r=s}^{n} V_r(\zeta) I_r, \qquad s \geq m+1,
\end{equation}
the function $H_s \colon Q_{u_s} \rightarrow \mathbb{C}$ defined as $H_s(\xi) := V_s(\zeta)$ for all
$\xi = \tau + i \eta=f_{u_s}(\zeta)$ for $\zeta\in \Sigma_{u_s}$ and $\tau,\eta\in\mathbb{R}$, where
the number $u_s$ is defined in the multiplication rule 3 for the basis $\{I_r\}_{r=1}^{n}$, has the partial derivatives $\partial H_s/\partial \tau$, $\partial H_s/\partial \eta$ at the point $\xi_0=f_{u_s}(\zeta_0)$, and equality \eqref{clas-CR-cond} holds for $u=s$.
}



\textbf{Proof.} (a) Let a vector $h_l$ be collinear to the respective vector $\zeta_l - \zeta_0$ for $l=1,2,3,4$, and
let $\theta_l := f_u(h_l) \in \{1, i, -1, -i\}$. 

Since the function $\Phi$ is differentiable in the sense of
G\^ateaux at a point $\zeta_0$, we have
\begin{equation} \label{Gprz-zeta0}
\lim\limits_{\delta\rightarrow 0^+} \frac{\Phi(\zeta_0+\delta
h_l)-\Phi(\zeta_0)}{\delta}= h_l\,\Phi_G'(\zeta_0).
\end{equation}

Applying the continuous multiplicative functional
$f_u$ to the both parts of equality \eqref{Gprz-zeta0}, we get
\[
\lim\limits_{\delta\rightarrow 0^+} \frac{V_u(\zeta_0+\delta h_l)-V_u(\zeta_0)}{\delta}= \theta_l\,f_u\big(\Phi_G'(\zeta_0)\big)
\]
that can be rewritten as
\begin{equation} \label{cor-Gprz-zeta0}
\lim\limits_{\delta\rightarrow 0^+} \frac{H_u(\xi_0+\delta \theta_l)-H_u(\xi_0)}{\delta}= \theta_l\,f_u\big(\Phi_G'(\zeta_0)\big).
\end{equation}

Equalities \eqref{cor-Gprz-zeta0} for $l=1,2,3,4$  yield the following equalities:
\[
\frac{\partial H_u}{\partial \tau}= \lim\limits_{\delta\rightarrow 0^+} \frac{H_u(\xi_0+\delta)-H_u(\xi_0)}{\delta} = f_u\big(\Phi_G'(\zeta_0)\big) = \lim\limits_{\delta\rightarrow 0^+} \frac{H_u(\xi_0-\delta)-H_u(\xi_0)}{-\delta},
\]
\[
\frac{\partial H_u}{\partial \eta} = \lim\limits_{\delta\rightarrow 0^+} \frac{H_u(\xi_0+\delta i)-H_u(\xi_0)}{\delta} = i\,f_u\big(\Phi_G'(\zeta_0)\big) = \lim\limits_{\delta\rightarrow 0^+} \frac{H_u(\xi_0-\delta i)-H_u(\xi_0)}{-\delta},
\]
which also imply the equality \eqref{clas-CR-cond}.


(b) Let $\Sigma_{u_s}\subset\Omega$ be a two-dimensional surface  containing four segments
$s[\zeta_{0}, \zeta_{1}]$, $s[\zeta_{0}, \zeta_{2}]$, $s[\zeta_{0}, \zeta_{3}]$, $s[\zeta_{0}, \zeta_{4}]$ such that their images
under the mapping $f_{u_s}$ are subsets of the segments $s[\xi_0,\xi_0 + 1]$, $s[\xi_0,\xi_0 + i]$, $s[\xi_0,\xi_0 - 1]$, $s[\xi_0,\xi_0 - i]$, respectively, where $\xi_0 = f_{u_s}(\zeta_0)$, and, in addition, let the restriction of the functional $f_{u_s}$ to the surface $\Sigma_u$ carry out a one-to-one mapping of this surface onto the domain $Q_{u_s}:=\{\xi=f_{u_s}(\zeta) : \zeta\in \Sigma_{u_s}\}$.

Let a vector $h_l$ be collinear to the respective vector $\zeta_l - \zeta_0$ for $l=1,2,3,4$, and
let $\theta_l := f_{u_s}(h_l) \in \{1, i, -1, -i\}$.


Since the function $\Phi$ is differentiable in the sense of G\^ateaux at a point $\zeta_0\in E_k$, we have equality \eqref{Gprz-zeta0} and, moreover,
\[ \Phi_G'(\zeta_0)=\frac{\partial\Phi}{\partial x_1}(\zeta_0)= \sum\limits_{r=s}^{n} \frac{\partial V_r}{\partial x_1}(\zeta_0)\,I_r\,.\]

Let us substitute the expression  $h_l=\theta_l I_{u_s}+\rho_l$ into equality \eqref{Gprz-zeta0}, where $\rho_l$ is some vector belonging to the ideal $\mathcal{I}_{u_s}$. Taking into account the multiplication rules for the Cartan basis and
the uniqueness of decomposition of elements of the algebra $\mathbb{A}_n^m$ with respect to the basis,
we equate complex-valued functions for $I_s$ on the both sides of equality \eqref{Gprz-zeta0} and get
\[
\lim\limits_{\delta\rightarrow 0^+} \frac{V_s(\zeta_0+\delta h_l)-V_s(\zeta_0)}{\delta}= \theta_l\,\frac{\partial V_r}{\partial x_1}(\zeta_0)\,.
\]
It can be rewritten as
\[\lim\limits_{\delta\rightarrow 0^+} \frac{H_s(\xi_0+\delta \theta_l)-H_s(\xi_0)}{\delta}= \theta_l\,\frac{\partial V_r}{\partial x_1}(\zeta_0).\]

Now, similarly to the part (a) of the proof, we can conclude that the partial derivatives $\partial H_s/\partial \eta$, $\partial H_s/\partial \tau$ exist at the point $\xi_0=f_{u_s}(\zeta_0)$ and equality \eqref{clas-CR-cond} holds for $u=s$.

Lemma is proved.

 \vskip2mm

\textbf{Lemma 2.}\label{lem-1} \textit{Suppose that the condition \eqref{cond-on-Ek} is satisfied and all intersections of a domain $\Omega\subset E_k$ with $(k-2)$-planes parallel to $M_u$ are linearly connected. Suppose also that 
 $\Phi:\Omega\rightarrow\mathbb{A}_n^m$ is a function locally bounded and differentiable in the sense of G\^ateaux in $\Omega$. If $\zeta_{1},\zeta_{2}\in\Omega$ and $\zeta_{2}-\zeta_{1}\in M_u$,
then
\[\Phi(\zeta_2)-\Phi(\zeta_1)\in\mathcal{I}_u\,.\]
Moreover, if the function $\Phi$ 
is of the form \eqref{Phi-expansion-rad}, then for $\zeta_{1},\zeta_{2}\in\Omega$ such that $\zeta_{2}-\zeta_{1}\in M_{u_s}$, the following inclusion holds: }
\begin{equation}\label{inclus-rad}
\Phi(\zeta_2)-\Phi(\zeta_1) \in \Biggr\{\sum\limits_{r=s+1}^n\lambda_rI_r:\lambda_r\in \mathbb{C}\Biggr\}\,.
\end{equation}

\textbf{Proof.} Let $\zeta_{1},\zeta_{2}\in\Omega$ and $\zeta_{2}-\zeta_{1}\in M_u$. Consider
the $(k-2)$-plane $L$ containing $\zeta_{1}$ and parallel to $M_u$. It is evident that  $\zeta_{2}\in L$.

Inasmuch as
$L \cap \Omega$ is linearly connected,
we can choose a curve $\gamma$ that lies in $L \cap \Omega$ and connects the points $\zeta_{1}$ and $\zeta_{2}$. By virtue of the compactness of the curve $\gamma$, there exists a finite set of open balls in $\Omega$ with centers lying on $\gamma$, which are covering this curve.

Let us add the points $\zeta_{1}$ and $\zeta_{2}$ to the centers of the mentioned balls and number them as
$\psi_1,\psi_2,\dots,\psi_m$
in the order in which they are encountered when moving along the curve $\gamma$ starting from the point $\psi_1=\zeta_1$ and ending with the point $\psi_m=\zeta_2$.

Consider the polygonal chain consisting of the line segments $s[\psi_j, \psi_{j+1}]$, $j=1,2,\dots,m-1$, each of which is parallel to some vector belonging to $M_u$ owing to the fact that $\psi_j, \psi_{j+1}\in\gamma\subset L$ and $L$ is parallel to $M_u$.





It is clear that there exists $\rho>0$ such that
for $j=1,2,\dots,m-1$ the convex neighborhood
\[\Omega_j:=\{\xi\in\Omega \, : \, \|\xi-\zeta\|<\rho, \,\,\, \zeta\in s[\psi_j, \psi_{j+1}]\,\} \]
of the line segment $s[\psi_j, \psi_{j+1}]$
is contained in the mentioned cover of the curve $\gamma$ by a finite set of open balls in $\Omega$.

Now, let us prove the relations $\Phi(\psi_{j+1})-\Phi(\psi_j)\in\mathcal{I}_u$
equivalent to the equalities $f_u(\Phi(\psi_{j+1})) = f_u(\Phi(\psi_{j}))$ for $j=1,2,\dots,m-1$.



Since 
the condition \eqref{cond-on-Ek} is satisfied, there exists an element $e_2^*\in E_k$
such that $f_u(e_2^*)=i$. Consider the lineal span $E_3^*:=\{\zeta=xe_1^*+ye_2^*+ze_3^* : x,y,z\in\mathbb{R}\}$ of the vectors $e_1^*=1$, $e_2^*$ and $e_3^*=\psi_{j+1}-\psi_j$.  
Introduce the convex neighborhood $\Omega_j^*:=\Omega_j\cap E_3^*$
of the line segment $s[\psi_j, \psi_{j+1}]$ in the three-dimensional real vector space $E_3^*$.

Let us construct in $\Omega_j^*$ two surfaces $\Upsilon_j$ and $\Sigma_j$
satisfying the following conditions (cf. \cite[p. 143]{Pla-Shpak-mono}):
\begin{itemize}
\item $\Upsilon_j$ and $\Sigma_j$ have the same edge;

\item the surface $\Upsilon_j$ contains the point $\psi_{j}$ and
the surface $\Sigma_j$ contains the point $\psi_{j+1}$;

\item restrictions of the functional $f_u$ onto the sets $\Upsilon_j$ and $\Sigma_j$ 
are one-to-one mappings of the mentioned sets onto the same domain $Q_j$ of the complex plane.

\item any line segment $s[\xi_0,\xi] \subset Q_j$ contains a line subsegment $s[\xi_0,\xi_1]$ such that its preimages on the surfaces $\Upsilon_j$ and $\Sigma_j$ over the mapping $f_u$ are the line segments as well.
\end{itemize}

As the surface $\Upsilon_j$, we can take an equilateral triangle having the center $\psi_{j}$ and, in addition, the plane of this triangle is perpendicular to the line segment $s[\psi_j, \psi_{j+1}]$ in $E_3^*$. As the surface $\Sigma_j$, we can take the lateral surface of the pyramid with the base $\Upsilon_j$ and the apex $\psi_{j+1}$.

Consider the function $V_u : \Omega\rightarrow \mathbb{C}$ defined by the equality $V_u(\zeta):= f_u \big(\Phi(\zeta)\big)$ for all $\zeta\in\Omega$. 

Now, the relation $V_u(\psi_{j+1})=V_u(\psi_{j})$ can be proved in a similar way as Lemma 5.3
 \cite{Sh-co} or Lemma 8.5 \cite{Pla-Shpak-mono}.

For each $\xi \in Q_j$, we define two complex-valued functions $H_1$ and $H_2$ as follows
$$H_1(\xi) := V_u(\zeta) \quad \mbox{for} \,\,\, \zeta \in \Upsilon_j,$$
$$H_2(\xi) := V_u(\zeta) \quad \mbox{for} \,\,\, \zeta \in \Sigma_j,$$
where the correspondence between the points $\zeta$ and $\xi \in Q_j$ is determined by the relation $f_u(\zeta) = \xi$.

Each of the functions $H_1, H_2$ satisfies the classical Cauchy--Riemann condition of the form \eqref{clas-CR-cond} at all points $\xi \in Q_j$
due to the part (a) of Lemma 1. Therefore, by virtue of Tolstoff’s result \cite{Tolstov}, the functions $H_1$ and $H_2$ are holomorphic in the domain $Q_j$.

Since, in addition, $H_1, H_2$ are continuous in the closure of domain $Q_j$ and $H_1(\xi) \equiv H_2(\xi)$ on the boundary of $Q_j$, this identity is also fulfilled everywhere in $Q_j$.
Therefore, $V_u(\psi_{j+1})=V_u(\psi_{j})$ and $\Phi(\psi_{j+1})-\Phi(\psi_j)\in\mathcal{I}_u$. Hence, taking into account the equality
\[\Phi(\zeta_2)-\Phi(\zeta_1)=\sum\limits_{j=1}^{m-1} \big(\Phi(\psi_{j+1})-\Phi(\psi_j)\big), \]
we get $\Phi(\zeta_2)-\Phi(\zeta_1)\in\mathcal{I}_u$.

In the case where the function $\Phi$ is of
the form \eqref{Phi-expansion-rad}, in the above reasoning,
we can consider the function $V_s$ instead of $V_u$, the set $M_{u_s}$ instead of $M_{u}$ and use the part (b) of Lemma 1. 
In such a way, for $\zeta_{1},\zeta_{2}\in\Omega$ such that $\zeta_{2}-\zeta_{1}\in M_{u_s}$, we get the equality
 $V_s(\zeta_2) - V_s(\zeta_1) = 0$ that yields inclusion \eqref{inclus-rad}.

  Lemma is proved.

 \vskip 2mm

Note that the condition of Lemma 2 on the linear connection of all intersections of the domain $\Omega\subset E_k$ with $(k-2)$-planes parallel to $M_u$ is a weakening of the convexity condition with respect to the set of directions $M_u$ in the case $k > 3$, which is accepted in the papers \cite{Pl-Pukh-Analele,Sh-co-Zb}. At the same time, these conditions are equivalent in the case $k = 3$.

 \vskip2mm

\subsection{A constructive description of locally bounded and differentiable in the sense of G\^ateaux functions}


A constructive description of monogenic functions given in a domain $\Omega\subset E_3$ is obtained in the paper \cite{Sh-co}
(see also Theorem 8.2 in \cite{Pla-Shpak-mono}) by means of holomorphic functions of complex variables.
This result is generalized in the paper \cite{Sh-co-Zb} to the case of a domain $\Omega\subset E_k$, where $2\leq k\leq 2n$.
As a consequence of such descriptions, it can be established that every function monogenic in $\Omega$
is differentiable in the sense of Lorch in the domain $\Omega$.

Now, we shall obtain the same constructive description of functions differentiable in the sense of G\^ateaux in a domain $\Omega\subset E_k$
provided that these functions are locally bounded, i.e., in comparison with the mentioned results of the papers \cite{Sh-co,Sh-co-Zb}, the condition of continuity of functions will be weakened.


 We introduce the linear operators $A_u$\,, $u=1,2,\ldots,m$, which
assign holomorphic functions $F_u:\,D_u\rightarrow\mathbb{C}$ to
every monogenic function
$\Phi:\Omega\rightarrow\mathbb{A}_n^m$ by the formula
\begin{equation}\label{def_op_A-k}
F_u(\xi_u):=f_u(\Phi(\zeta)),
\end{equation}
where $\xi_u=f_u(\zeta)$ 
and $\zeta\in\Omega$. It follows from Lemma 2
that the value $F_u(\xi_u)$ does
not depend on a choice of a point $\zeta$ for which
$f_u(\zeta)=\xi_u$.

 The following theorem is proved in a similar way as Theorem 1 \cite{Sh-co-Zb} (see also Theorem 8.2 in \cite{Pla-Shpak-mono}).

\vskip2mm

\textbf{Theorem 2.}\label{teo_1_konstruct_opys_A_n_m-k} \textit{Suppose that condition \eqref{cond-on-Ek} is satisfied and
all intersections of a domain $\Omega\subset E_k$ with planes parallel to $M_u$ are linearly connected.
Then every function $\Phi:\Omega\rightarrow\mathbb{A}_n^m$ locally bounded and
differentiable in the sense of G\^ateaux
can be expressed in the form
 \begin{multline}\label{Teor--1-k} 
\Phi(\zeta)=\sum\limits_{u=1}^mI_u\,\frac{1}{2\pi
i}\int\limits_{\Gamma_u} F_u(t)(te_1-\zeta)^{-1}\,dt\\
+\sum\limits_{s=m+1}^nI_s\,\frac{1}{2\pi i}\int\limits_
{\Gamma_{u_s}}G_s(t)(te_1-\zeta)^{-1}\,dt,
 \end{multline} 
where $F_u$ and $G_s$ are certain holomorphic functions in the
domains $D_u$ and $D_{u_s}$, respectively, and $\Gamma_u$ is a
closed rectifiable Jordan curve in $D_u$, which embraces the point
$\xi_u$ and contains no points $\xi_q$ for $q=1,2,\ldots,
m$, $q\neq u$.}
 \vskip2mm

\textbf{Proof.}
We define the holomorphic functions\, $F_u$\,, $u=1,2,\ldots,m$\,, by equality \eqref{def_op_A-k}.

Let us show that the values of the function
\begin{equation}\label{teor__2-k}
\Phi_0(\zeta):=\Phi(\zeta)-\sum\limits_{u=1}^mI_u\,\frac{1}{2\pi
i} \int\limits_{\Gamma_u}F_u(t)(te_1-\zeta)^{-1}\,dt
\end{equation}
belong to the radical $\mathcal{R}$,
i.e., $\Phi_0(\zeta)\in\mathcal{R}$ for all
$\zeta\in\Omega$.

As a consequence of equality (\ref{lem-rez}), we obtain the
equality
\begin{multline*}
I_u\,\frac{1}{2\pi i}
\int\limits_{\Gamma_u}F_u(t)(te_1-\zeta)^{-1}\,dt=I_u\,\frac{1}{2\pi
i} \int\limits_{\Gamma_u}\frac{F_u(t)}{t-\xi_u}\,dt+\\
 +\frac{1}{2\pi i}\sum\limits_{s=m+1}^{n}
 \sum\limits_{r=2}^{s-m+1}\int\limits_{\Gamma_u}\frac{F_u(t)Q_{r,s}}
 {\left(t-\xi_{u_{s}}\right)^r}\,dt
 \,I_{s}\,I_u\,,
\end{multline*}
which implies the equality
\begin{equation}\label{teor__3-k}
f_u\left(\sum\limits_{u=1}^mI_u\,\frac{1}{2\pi i}
\int\limits_{\Gamma_u}F_u(t)(te_1-\zeta)^{-1}\,dt\right)=F_u(\xi_u).
\end{equation}

Acting on equality (\ref{teor__2-k}) by the functional $f_u$
and taking into account relations (\ref{def_op_A-k}) and (\ref{teor__3-k}),
we get the equality
$$f_u(\Phi_0(\zeta))=F_u(\xi_u)-F_u(\xi_u)=0
$$
for all $u=1,2,\ldots,m$, i.e., $\Phi_0(\zeta)\in\mathcal{R}$.

Therefore, the function $\Phi_0$ is of the form
\[\Phi_{0}(\zeta)=\sum\limits_{s=m+1}^{n} V_{s}(\zeta)\,I_s\,,\]
where $V_{s}:\Omega\rightarrow\mathbb{C}$.


It follows from Lemma 2 that
$V_{m+1}(\zeta)\equiv G_{m+1}(\xi_{u_{m+1}})$ for all $\zeta\in\Omega$, where 
$\xi_{u_{m+1}}=f_{u_{m+1}}(\zeta)$ and $G_{m+1} : D_{u_{m+1}}\rightarrow\mathbb{C}$.
In addition, the function $G_{m+1}$ satisfies the classical Cauchy--Riemann condition of the form \eqref{clas-CR-cond} at all points $\xi_{u_{m+1}}\in D_{u_{m+1}}$ due to the part (b) of Lemma 1. 
Hence, by virtue of Tolstoff’s result \cite{Tolstov}, the function $G_{m+1}$ is holomorphic in the domain $D_{u_{m+1}}$\,.
Therefore,
\begin{equation}\label{teor__4}
\Phi_0(\zeta)=G_{m+1}(\xi_{u_{m+1}})\,I_{m+1}+
\sum\limits_{s=m+2}^{n} V_{s}(\zeta)\,I_s\,.
\end{equation}

By virtue of expansion (\ref{lem-rez}), we have the representation
\begin{equation}\label{teor__5}
I_{m+1}\,\frac{1}{2\pi
i}\int\limits_{\Gamma_{u_{m+1}}}G_{m+1}(t)(te_1-\zeta)^{-1}\,dt=
G_{m+1}(\xi_{u_{m+1}})\,I_{m+1}+\Psi(\zeta),
\end{equation}
where $\Psi(\zeta)$ is a function taking values in the set
\[\bigg\{\sum_{s=m+2}^n\lambda_s\,I_s:\lambda_s\in\mathbb{C}\bigg\}.\]

Now, consider the function
 $$\Phi_1(\zeta):=\Phi_0(\zeta)-I_{m+1}\,\frac{1}{2\pi
i}\int\limits_{\Gamma_{u_{m+1}}}G_{m+1}(t)(te_1-\zeta)^{-1}\,dt.$$
In view of the relations (\ref{teor__4}), (\ref{teor__5}),
the function $\Phi_1$ can be represented as
$$\Phi_{1}(\zeta)=\sum\limits_{s=m+2}^n \widetilde
V_{s}(\zeta)\,I_s\,,$$
where $\widetilde V_{s}:\Omega\rightarrow\mathbb{C}$\,.

Inasmuch as the function $\Phi_1$ is locally bounded and
differentiable in the sense of G\^ateaux in $\Omega$, the function $\widetilde V_{m+2}$ is
similar to the function $V_{m+1}$, and we can state that
$\widetilde V_{m+2}(\zeta)\equiv G_{m+2}(\xi_{u_{m+2}})$ for all $\zeta\in\Omega$, where
$\xi_{u_{m+2}}=f_{u_{m+2}}(\zeta)$ and $G_{m+2} : D_{u_{m+2}}\rightarrow\mathbb{C}$ is a function holomorphic in the domain
$D_{u_{m+2}}$\,.

In such a way, step by step, considering the functions
$$\Phi_j(\zeta):=\Phi_{j-1}(\zeta)-I_{m+j}\,\frac{1}{2\pi
i}\int\limits_{\Gamma_{u_{m+j}}}G_{m+j}(t)(te_1-\zeta)^{-1}\,dt$$
 for $j=2,3,\dots,n-m-1$, we get representation
(\ref{Teor--1-k}) of the function $\Phi$.

The theorem is proved.

 \vskip2mm

As a result of Theorem 2, 
we get the statement of Theorem 1. Indeed, it follows from representation \eqref{Teor--1-k} that the function $\Phi$ is
differentiable in the sense of Lorch in the domain $\Omega$, which in turn implies that the function $\Phi$ is monogenic
in $\Omega$.

 \vskip2mm

\textbf{Acknowledgments.} 
This work was supported by a grant from the Simons Foundation (1290607, S.A.P., V.S.Sh., M.V.T.).

\end{document}